\documentclass[11pt]{article}
\usepackage{amsmath, amssymb, amscd, epsfig}

\newtheorem{thm}{Theorem}

\newtheorem{lm}[thm]{Lemma}

\newtheorem{rem}[thm]{Remark}

\renewcommand{\phi}{\varphi}
\renewcommand{\epsilon}{\varepsilon}

\newcommand{\BB}{\mathbb}
\newcommand{\g}{\mathfrak}

\newcommand{\separate}{\vskip5pt}

\begin{document}

\title{\bf A Localization Argument for Characters of Reductive Lie Groups:
An Introduction and Examples}
\author{Matvei Libine}
\maketitle

\begin{abstract}
In this article I describe my recent geometric localization argument
dealing with actions of NONcompact groups which
provides a geometric bridge between two
entirely different character formulas for reductive Lie
groups and answers the question posed in \cite{Sch}.

A corresponding problem in the compact group setting was
solved by N.~Berline, E.~Getzler and M.~Vergne in \cite{BGV}
by an application of the theory of equivariant forms
and, particularly, the fixed point integral localization formula.

This localization argument seems to be the first successful attempt
in the direction of building a similar theory for integrals of differential
forms, equivariant with respect to actions of noncompact groups.

I will explain how the argument works in the $SL(2,\BB R)$ case,
where the key ideas are not obstructed by technical details and
where it becomes clear how it extends to the general case.
The general argument appears in \cite{L}.

I put every effort to make this article broadly accessible.
Also, although this article mentions characteristic cycles of sheaves,
I do not assume that the reader is familiar with this notion.
\end{abstract}

\tableofcontents

\begin{section}
{Introduction}
\end{section}

For motivation, let us start with the case of a compact group.
Thus we consider a connected compact group $K$ and a
maximal torus $T \subset K$.
Let $\g k_{\BB R}$ and $\g t_{\BB R}$ denote the Lie algebras
of $K$ and $T$ respectively, and $\g k$, $\g t$ be their
complexified Lie algebras.
Let $\pi$ be a finite-dimensional representation of $K$,
that is $\pi$ is a group homomorphism $K \to Aut(V)$.
One of the most important invariants associated to a
representation is its {\em character}. It is a function
on $K$ defined by
$\Theta_{\pi}(x) =_{\text{def}} \operatorname{tr} (\pi(x))$,
$x \in K$.
Any finite-dimensional representation is completely determined
(up to isomorphism) by its character.
Recall the exponential map $\exp: \g k_{\BB R} \to K$.
If $K$ is a matrix subgroup of some $GL(n)$, then
$\exp(A) = e^A$. We use this exponential map to define the
{\em character on the Lie algebra} of the representation
$(\pi,V)$:
$$
\theta_{\pi}= (\det \exp_*)^{1/2} \exp^* \Theta_{\pi}.
$$
It is a smooth bounded function on $\g k_{\BB R}$.
Because $K$ is connected and compact, the exponential map is
surjective and generically non-singular.
Thus $\theta_{\pi}$ still determines the representation.

Now let us assume that the representation $\pi$ is irreducible,
that is there are no proper $K$-invariant subspaces of $V$.
One reason for introducing this factor $(\det \exp_*)^{1/2}$
in the definition of $\theta_{\pi}$ is that for irreducible
representations $\pi$ the character $\theta_{\pi}$ becomes
an eigenfunction of the algebra of conjugation-invariant
constant coefficient differential operators on
the Lie algebra $\g k_{\BB R}$.

There are two entirely different character formulas
for $\theta_{\pi}$ -- the {\em Weyl character formula}
and {\em Kirillov's character formula}.
Recall that the irreducible representations of $K$
can be enumerated by their {\em highest weights} which are
elements of the {\em weight lattice} $\Lambda$ in $i \g t_{\BB R}^*$
intersected with a chosen {\em Weyl chamber}.
Let $\lambda=\lambda(\pi) \in i \g t_{\BB R}^*$
denote the highest weight corresponding to $\pi$.
The group $K$ acts on its own Lie algebra
$\g k_{\BB R}$ by adjoint representation $Ad$.
Let $W=N_K(\g t_{\BB R})/T$,
where $N_K(\g t_{\BB R})$ is the normalizer
of $\g t_{\BB R}$ in $K$.
The set $W$ is a finite group called the {\em Weyl group};
it acts on $\g t_{\BB R}$ and hence on $i \g t_{\BB R}^*$.
We can choose a positive definite inner product
$\langle \cdot,\cdot\rangle$ on $i \g t_{\BB R}^*$ invariant under $W$.
Then the Weyl character formula can be stated as follows:
$$
\theta_{\pi}|_{\g t_{\BB R}}(t)=\sum_{w \in W} \frac {e^{w(\lambda+\rho)(t)}}
{\prod_{\alpha \in \Phi, \langle w(\lambda+\rho),\alpha\rangle>0} \alpha(t)},
$$
where $\Phi \subset i \g t_{\BB R}^*$ is the {\em root system}
of $\g k_{\BB R}$ and
$\rho \in i \g t_{\BB R}^*$ is a certain small vector independent
of $\pi$. Because $\theta_{\pi}$ is $Ad(K)$-invariant and
every $Ad(K)$-orbit in $\g k_{\BB R}$ meets $\g t_{\BB R}$,
this formula completely determines $\theta_{\pi}$.

Kirillov's character formula provides a totally different expression
for the irreducible characters on $\g k_{\BB R}$.
The splitting
$\g k_{\BB R}=\g t_{\BB R} \oplus [\g t_{\BB R},\g k_{\BB R}]$
(Cartan algebra $\oplus$ root spaces)
induces a dual splitting of the vector space $i \g k_{\BB R}^*$,
which allows us to think of $\lambda$ and $\rho$ as lying in
$i \g k_{\BB R}^*$. The adjoint action of $K$ on
$\g k_{\BB R}$ has a dual action on $i \g k_{\BB R}^*$
called coadjoint representation. We define
$$
\Omega_{\lambda+\rho} = \text{ $K$-orbit of $\lambda+\rho$
in $i \g k_{\BB R}^*$}.
$$
It will be convenient to define the Fourier transform $\hat \phi$
of a test function $\phi \in {\cal C}^{\infty}_c (\g k_{\BB R})$
without the customary factor of $i=\sqrt{-1}$ in the exponent,
as a function on $i \g k_{\BB R}^*$:
$$
\hat \phi(\zeta) = \int_{\g k_{\BB R}}
\phi(x)e^{\langle \zeta, x \rangle} dx.
$$
Then Kirillov's character formula describes $\theta_{\pi}$ as a
distribution on $\g k_{\BB R}$:
$$
\int_{\g k_{\BB R}} \theta_{\pi} \phi dx =
\int_{\Omega_{\lambda+\rho}} \hat \phi d\beta,
$$
where $d\beta$ is the measure induced by the canonical symplectic
structure of $\Omega_{\lambda+\rho}$.
In other words,
$$
\hat \theta_{\pi} = \text{ integration over $\Omega_{\lambda+\rho}$}.
$$
Kirillov calls this the ``universal formula'' for irreducible characters.

\separate

The geometric relationship between these two formulas is even
more striking. As a homogeneous space, $\Omega_{\lambda+\rho}$
is isomorphic to the {\em flag variety} $X$, i.e the variety
of Borel subalgebras
$\g b \subset \g k = \g k_{\BB R} \otimes_{\BB R} \BB C$.
The space $X$ is a smooth complex projective variety which is also
isomorphic to $K/T$ as a homogeneous space.
The Borel-Weil-Bott theorem can be regarded as an explicit
construction of a holomorphic $K$-equivariant line
bundle ${\cal L}_{\lambda} \to X$ such that the resulting
representation of $K$ in the cohomology groups is:
\begin{eqnarray*}
&& H^p(X, {\cal O}({\cal L}_{\lambda}))=0 \text{\quad if $p \ne 0$},  \\
&& H^0(X, {\cal O}({\cal L}_{\lambda})) \simeq \pi.
\end{eqnarray*}
Then the Weyl character formula is a consequence of
the Atiyah-Bott fixed point formula or the Lefschetz fixed point formula.
On the other hand, N.~Berline, E.~Getzler
and M.~Vergne proved in \cite{BGV}
Kirillov's character formula using the integral
localization formula for $K$-equivariant forms.
They showed that the right hand side of Kirillov's character formula
equals the right hand side of the Weyl character formula.

\separate

Equivariant forms were introduced in 1950 by Henri Cartan.
There are many good texts on this subject including
\cite{BGV} and \cite{GS}.

Let $K$ act on a smooth manifold $M$, and let $\Omega^{\bullet}(M)$
denote the algebra of smooth differential forms on $M$.
A $K$- equivariant form is a smooth map
$$
\omega: \g k_{\BB R} \to \Omega^{\bullet}(M)
$$
whose image need not lie entirely in any single degree component of 
$\Omega^{\bullet}(M)$ and which is $K$-invariant,
i.e. for all elements $k \in \g k_{\BB R}$ and $\tilde k \in K$ we have
$$
\omega(k) = (\tilde k)^{-1} \cdot \omega \bigl( Ad(\tilde k) k \bigr).
$$

We say that an equivariant form $\omega$ is {\em equivariantly closed}
if $d_{equiv}(\omega)=0$, where
$$
\bigl( d_{equiv}(\omega) \bigr) (k) = d(\omega(k)) + \iota (k) \omega(k).
$$
Here, the first term is the ordinary deRham differential of $\omega(k)$
and the last term denotes the contraction of the differential form
$\omega(k)$ by the vector field on $M$ generated by the infinitesimal
action of $k$.

If $N \subset M$ is a submanifold, then
$$
\int_N \omega =_{\text{def}}
\int_N \text{component of $\omega$ of degree $\dim N$};
$$
it is a function on $\g k_{\BB R}$. Then the localization formula
reduces integration of an equivariantly closed form to
summation over the zeroes of the vector field in $M$ generated
by $k \in \g k_{\BB R}$.
In other words, it expresses a global object such as integral of a
differential form in terms of quantities which can be calculated
locally at the zeroes of the vector field.
It is crucial for the localization formula to hold that the group $K$
is compact.

\separate

Now let $G_{\BB R}$ be a connected, linear, reductive Lie group.
We let $\g g_{\BB R}$ denote its Lie algebra.
Then most representations of interest have infinite dimension.
We always consider representations on complete, locally convex
Hausdorff topological vector spaces and require that the action
of $G_{\BB R}$ is continuous.
Let $K$ be a maximal compact subgroup of $G_{\BB R}$.
A reasonable category of representations consists of
{\em admissible} representations of {\em finite length}.
(A representation $\pi$ has finite length if every increasing chain
of closed, invariant subspaces breaks off after finitely many steps;
$\pi$ is admissible if its restriction to $K$ contains any irreducible
representation of $K$ at most finitely often.) Admissibility
is automatic for irreducible unitary representations.
Although trace of a linear operator in an infinite-dimensional
space cannot be defined in general, it is still possible to
define a character $\theta_{\pi}$ as an
$Ad(G_{\BB R})$-invariant distribution on $\g g_{\BB R}$.

M.~Kashiwara and W.~Schmid in their paper \cite{KSch}
generalize the Borel-Weil-Bott construction.
(Recall that $X$ denotes the flag variety of Borel subalgebras of the
complexified Lie algebra $\g g = \g g_{\BB R} \otimes_{\BB R} \BB C$.)
Instead of line bundles on $X$ they consider
$G_{\BB R}$-equivariant sheaves ${\cal F}$ and, for each integer
$p \in \BB Z$, they define representations of $G_{\BB R}$ in
$\operatorname{Ext}^p({\cal F},{\cal O})$.
Such representations turn out to be admissible of finite length.
Then W.~Schmid and K.~Vilonen prove in \cite{SchV2}
two character formulas for these representations --
the fixed point character formula and the integral character formula.
In the case when $G_{\BB R}$ is compact, the former reduces to
the Weyl character formula and the latter -- to Kirillov's character
formula. The fixed point formula was conjectured by M.~Kashiwara
in \cite{K}, and its proof uses a
generalization of the Lefschetz fixed point formula to sheaf
cohomology due to M.~Goresky and R.~MacPherson in \cite{GM}.
On the other hand, W.~Rossmann in \cite{R} established existence
of an integral character formula over an unspecified
Borel-Moore cycle.
W.~Schmid and K.~Vilonen prove the integral character formula
where integration takes place over the {\em characteristic cycle} of
${\cal F}$, $Ch({\cal F})$, and their proof
depends totally on representation theory.

$Ch({\cal F})$ is a conic Lagrangian cycle in the cotangent space
$T^*X$ associated to the sheaf ${\cal F}$.
Characteristic cycles were introduced by M.~Kashiwara and their
definition can be found in \cite{KaScha}. On the other hand,
W.~Schmid and K.~Vilonen give a geometric way to understand
characteristic cycles in \cite{SchV1}.

In this article we only deal with characteristic cycles when
considering the important cases in the situation of $SL(2,\BB R)$.
Each time a characteristic cycle appears here it is described
explicitly and, therefore, the reader is not expected to be
familiar with this notion.

\separate

The equivalence of these two formulas can be stated in terms
of the sheaves ${\cal F}$ alone, without any reference to
their representation-theoretic significance.
In the announcement \cite{Sch} W.~Schmid posed a question:
``Can this equivalence be seen directly without a detour to
representation theory, just as in the compact case.''

In this article I provide such a geometric link.
I introduce a localization technique which allows
to localize integrals to the zeroes of vector fields on $X$
generated by the infinitesimal action of $\g g_{\BB R}$.
Thus, in addition to a representation-theoretical result,
we obtain a whole family of examples where it is possible
to localize integrals to fixed points with respect to an
action of a noncompact group.
As far as I am aware, this localization argument is the first
successful attempt in the direction of building an integral
localization theory for noncompact groups.

\separate

The following convention will be in force throughout these notes:
whenever $A$ is a subset of $B$,
we will denote the inclusion map $A \hookrightarrow B$ by
$j_{A \hookrightarrow B}$.

\separate

\begin{section}
{Two Character Formulas}
\end{section}

In these notes we try to keep the same notations as W.~Schmid and K.~Vilonen
use in \cite{SchV2} as much as possible.
That is they fix a connected, complex algebraic,
reductive group $G$ which is defined over $\BB R$.
The representations they consider are representations of a real form
$G_{\BB R}$ of $G$ -- in other words, $G_{\BB R}$ is a subgroup of
$G$ lying between the group of real points $G(\BB R)$ and the
identity component $G(\BB R)^0$.
They regard $G_{\BB R}$ as a reductive Lie group and
denote by $\g g$ and $\g g_{\BB R}$ the Lie algebras of
$G$ and $G_{\BB R}$ respectively,
they also denote by $X$ the flag variety of $G$.

If $g \in \g g$ is an element of the Lie algebra,
we denote by $\operatorname{VF}_g$ the vector field on $X$
generated by $g$: if $x \in X$ and $f \in {\cal C}^{\infty}(X)$, then
$$
\operatorname{VF}_g (x)f=
\frac d{d\epsilon} f(\exp(\epsilon g) \cdot x)|_{\epsilon=0}.
$$
We call a point $x \in X$ a {\em fixed point of $g$} if the vector field
$\operatorname{VF}_g$ on $X$ vanishes at $x$, i.e.
$\operatorname{VF}_g (x)= 0$.

In this paragraph we explain the general picture, but since objects
mentioned here will not play any role in what follows they will not
be defined, rather the reader is referred to \cite{SchV2}.
W.~Schmid and K.~Vilonen denote by $\g h$ the universal Cartan
algebra. They pick an element $\lambda \in \g h^*$ and introduce
the ``$G_{\BB R}$-equivariant derived category on $X$ with twist
$(\lambda - \rho)$'' denoted by
$\operatorname{D}_{G_{\BB R}}(X)_{\lambda}$.
They also introduce ${\cal O}_X(\lambda)$,
the twisted sheaf of holomorphic functions on $X$,
with twist $(\lambda - \rho)$.
Then, for ${\cal F} \in \operatorname{D}_{G_{\BB R}}(X)_{-\lambda}$,
they define a virtual representation of $G_{\BB R}$
$$
\sum_p (-1)^p \operatorname{Ext}^p (\BB D {\cal F}, {\cal O}_X(\lambda)),
$$
where $\BB D {\cal F} \in \operatorname{D}_{G_{\BB R}}(X)_{\lambda}$
denotes the Verdier dual of ${\cal F}$.
It was shown in \cite{KSch} that each
$\operatorname{Ext}^p (\BB D {\cal F}, {\cal O}_X(\lambda))$
is admissible of finite length.
There are two formulas expressing the character $\theta$ of this
virtual representation as a distribution on $\g g_{\BB R}$.
We think of a character as a linear functional defined on the space of
smooth compactly supported differential forms $\phi$ on $\g g_{\BB R}$
of top degree, and write $\int_{\g g_{\BB R}} \theta \phi$ for the value
of $\theta$ at $\phi$.
In my article \cite{L} I start with the right hand side of
the integral character formula
\begin{equation}  \label{intformula}
\int_{\g g_{\BB R}} \theta \phi =
\frac 1{(2\pi i)^nn!} \int_{Ch({\cal F})}
\mu_{\lambda}^* \hat \phi (-\sigma+\pi^* \tau_{\lambda})^n
\end{equation}
and show that it is equivalent to the right hand side of the
fixed point character formula
\begin{equation}  \label{fpformula}
\int_{\g g_{\BB R}} \theta \phi = \int_{\g g_{\BB R}} F_{\theta} \phi,
\qquad
F_{\theta} (g) =
\sum_{k=1}^{|W|} 
\frac {m_{x_k(g)} e^{\langle g,\lambda_{x_k(g)} \rangle}}
{\alpha_{x_k(g),1}(g) \dots \alpha_{x_k(g),n}(g)},
\end{equation}
where $F_{\theta}$ is a function in $L^1_{loc} (\g g_{\BB R})$,
$x_1, \dots, x_n$ are the fixed points of $g$ and
integers $m_{x_k(g)}$'s are the local invariants of ${\cal F}$.

\separate

Because both character formulas depend on
${\cal F} \in \operatorname{D}_{G_{\BB R}}(X)_{-\lambda}$ only through
its characteristic cycle $Ch({\cal F})$,
as far as characteristic cycles concern, we can simply
replace ${\cal F}$ with a $G_{\BB R}$-equivariant sheaf on the flag
variety $X$ with the same characteristic cycle.
We will use the same notation ${\cal F}$ to denote this
$G_{\BB R}$-equivariant sheaf on $X$.
Let $n = \dim_{\BB C} X$, let $\pi: T^*X \twoheadrightarrow X$ be
the projection map, and equip $\g g_{\BB R}$ with some orientation.
Then $Ch({\cal F})$ is a $2n$-cycle in $T^*X$ which has real dimension $4n$.

We will make an elementary calculation of the integral
$$
\frac 1{(2\pi i)^nn!} \int_{Ch({\cal F})}
\mu_{\lambda}^* \hat \phi (-\sigma+\pi^* \tau_{\lambda})^n
$$
where $\phi$ is a smooth compactly supported differential
form on $\g g_{\BB R}$ of top degree,
$$
\hat \phi(\zeta) = \int_{\g g_{\BB R}} e^{\langle g, \zeta \rangle} \phi
\qquad (g \in \g g_{\BB R},\,\zeta \in \g g^*)
$$
is its Fourier transform
(without the customary factor of $i=\sqrt{-1}$ in the exponent),
$\mu_{\lambda}: T^*X \to \g g^*$
is the twisted moment map defined in \cite{SchV1}
and $\tau_{\lambda}$, $\sigma$ are 2-forms on
$X$ and $T^*X$ respectively defined in \cite{SchV2}.
$\sigma$ is the complex algebraic symplectic form on $T^*X$.
On the other hand, the precise definition of the form $\tau_{\lambda}$
will not be important. What will be important, however,
is that, for each $g \in \g g$, the $2n$-form on $T^*X$
\begin{equation}  \label{integrand}
e^{\langle g, \mu_{\lambda}(\zeta) \rangle}
(-\sigma+\pi^* \tau_{\lambda})^n
\end{equation}
is closed.

\separate

If $\lambda =0$, the twisted moment map $\mu_{\lambda}$ reduces to
the ordinary moment map $\mu_0=\mu$ defined by
$$
\mu(\zeta): g \mapsto \langle \zeta, \operatorname{VF}_g \rangle,
\qquad \zeta \in T^*X.
$$
The moment map $\mu$ takes values in the nilpotent cone
${\cal N^*} \subset \g g^*$.
In general, $\mu_{\lambda} = \mu + \lambda_x$, where
$\lambda_x$ is a function on $X$ mapping $x \in X$ into
$\lambda_x \in \g g^*$.

Suppose that $\lambda \in \g h^*$ is regular. For $\g g = \g{sl}(2,\BB C)$
this simply means $\lambda \ne 0$. Then $\mu_{\lambda}$ is a real analytic
diffeomorphism of $T^*X$ onto $\Omega_{\lambda} \subset \g g^*$
-- the orbit of $\lambda$ under the coadjoint action of $G$ on $\g g^*$.
Let $\sigma_{\lambda}$ denote the canonical $G$-invariant complex algebraic
symplectic form on $\Omega_{\lambda}$. Then
$e^{\langle g, \zeta \rangle} (\sigma_{\lambda})^n$ is a holomorphic
$2n$-form of maximal possible degree, hence closed.
It turns out that
$\mu_{\lambda}^* (\sigma_{\lambda}) = -\sigma+\pi^* \tau_{\lambda}$.
This shows that, for $\lambda$ regular, $g \in \g g$,
$$
e^{\langle g, \mu_{\lambda}(\zeta) \rangle}
(-\sigma+\pi^* \tau_{\lambda})^n =
\mu_{\lambda}^* \bigl( e^{\langle g, \zeta \rangle} (\sigma_{\lambda})^n \bigr)
$$
is a closed $2n$-form on $T^*X$.
Note that neither map $\mu_{\lambda}$ nor the form (\ref{integrand})
is holomorphic. Because the form
$e^{\langle g, \mu_{\lambda}(\zeta) \rangle}
(-\sigma+\pi^* \tau_{\lambda})^n$
depends on $\lambda$ real analytically and the set of regular elements
is dense in $\g h^*$, we conclude that the form in the equation
(\ref{integrand}) is closed.

\separate

\begin{section}
{Localization Argument for $SL(2,\BB R)$}
\end{section}

In this section we will deal with characteristic cycles of sheaves.
I remind it once again that all characteristic cycles which
appear here will be described explicitly and the reader is not expected
to be familiar with this notion.

\separate

There are three major obstacles to any geometric proof of the equivalence
of formulas (\ref{intformula}) and (\ref{fpformula}):
first of all, the characteristic cycle $Ch({\cal F})$ need not be smooth,
in fact, it may be extremely singular;
secondly, the integrand in (\ref{intformula}) is an equivariant form with
respect to some compact real form $U_{\BB R} \subset G$,
but $U_{\BB R}$ does not preserve $Ch({\cal F})$
(unless $Ch({\cal F})$ is a multiple of the zero section of $T^*X$
equipped with some orientation);
and, lastly, $Ch({\cal F})$ is not compactly supported.

\separate

In this section we will consider two different examples when
$G_{\BB R} = SL(2,\BB R)$.

Recall that when $G_{\BB R} = SL(2,\BB R)$, the flag variety $X$ is
just a 2-sphere $\BB C P^1$ on which $SL(2,\BB R)$ acts by
projective transformations.
There are exactly three $SL(2,\BB R)$-orbits: two open hemispheres
and one circle which is their common boundary.
Also, $n = \dim_{\BB C} X =1$.

\separate

To deal with the problem of noncompactness of $Ch({\cal F})$
we fix a norm $\|.\|$ on the space of linear functionals on
$\g{sl}(2,\BB R)$, $\g{sl}(2,\BB R)^*$. Then the moment map $\mu$
induces a vector bundle norm on $T^*X$: for $\zeta \in T^*X$
its norm will be $\|\mu(\zeta)\|$. We will use the same notation
$\|.\|$ for this norm too.

Let $\g{sl}(2,\BB R)'$ denote the set of regular semisimple elements
in $\g{sl}(2,\BB R)$: $\g{sl}(2,\BB R)'$ consists of those elements
in $\g{sl}(2,\BB R)$ which are diagonalizable over $\BB C$,
with distinct eigenvalues.

Since the complement of $\g{sl}(2,\BB R)'$ in $\g{sl}(2,\BB R)$
has measure zero, we can replace integration over $\g{sl}(2,\BB R)$
by integration over $\g{sl}(2,\BB R)'$. Then
\begin{multline}  \label{R}
\frac 1{2\pi i} \int_{Ch({\cal F})}
\mu_{\lambda}^* \hat \phi (-\sigma+\pi^* \tau_{\lambda}) \\
=\lim_{R \to \infty}
\frac 1{2\pi i}
\int_{\g{sl}(2,\BB R) \times (Ch({\cal F}) \cap \{\|\zeta\| \le R\})}
e^{\langle g, \mu_{\lambda}(\zeta) \rangle} \phi
(-\sigma+\pi^* \tau_{\lambda}) \\
=\lim_{R \to \infty} \frac 1{2\pi i}
\int_{\g{sl}(2,\BB R)' \times (Ch({\cal F}) \cap \{\|\zeta\| \le R\})}
e^{\langle g, \mu_{\lambda}(\zeta) \rangle} \phi
(-\sigma+\pi^* \tau_{\lambda}).
\end{multline}
(Of course, the orientation on
$\g{sl}(2,\BB R)' \times (Ch({\cal F}) \cap \{\|\zeta\| \le R\})$
is induced by the product orientation on
$\g{sl}(2,\BB R) \times Ch({\cal F})$.)

We will interchange the order of integration:
integrate over the characteristic cycle first and only then
perform integration over $\g{sl}(2,\BB R)'$.
As we already mentioned, the integrand
is an equivariant form with respect to some
compact real form $U_{\BB R} \subset G = SL(2,\BB C)$.
But $U_{\BB R}$ does not preserve $Ch({\cal F})$
(unless $Ch({\cal F})$ is a multiple of the zero section of $T^* \BB CP^1$
equipped with some orientation).
We regard the integral as an integral of a closed differential form
$e^{\langle g, \mu_{\lambda}(\zeta) \rangle}
\phi(-\sigma+\pi^* \tau_{\lambda})$
over a chain in
$\g{sl}(2,\BB R)' \times (T^*\BB CP^1 \cap \{\|\zeta\| \le R\})$.

Each regular semisimple element $g \in \g{sl}(2,\BB R)'$
has exactly two fixed points on $X=\BB CP^1$.
We will use the open embedding theorem of W.~Schmid and K.~Vilonen
(\cite{SchV1}) to construct a deformation of
$Ch({\cal F})$ into a simple cycle of the following kind:
$$
m_1 T^*_{x_1}X + m_2 T^*_{x_2}X,
$$
where $m_1,m_2$ are some integers,
$x_1,x_2$ are the points in $\BB CP^1$ fixed by $g$,
and each cotangent space $T^*_{x_k}X$ is given some orientation.
To ensure that the integral behaves well, we will stay during the
process of deformation inside the set
$$
\{(g,\zeta) \in \g{sl}(2,\BB R)' \times T^* \BB CP^1;\,
Re( \langle g,\mu(\zeta) \rangle) \le 0\}.
$$

\separate

The most important case is the split Cartan case,
i.e. the case when the eigenvalues of $g$ are real and
the unique Cartan algebra containing $g$ is split.
The compact Cartan case also needs to be considered, but
because it can be treated using the classical theory of
equivariant forms it is not so interesting.

In the split Cartan case one fixed point in $\BB CP^1$ is stable
(the vector field $\operatorname{VF}_g$ on $\BB CP^1$ points towards it)
and the other fixed point in $\BB CP^1$ is unstable (the vector field 
$\operatorname{VF}_g$ points away from it).
We can position $X=\BB CP^1$ in space so that the unstable fixed point
is the north pole $N$ and the stable fixed point is the south pole $S$.
Furthermore, we can position it so that the eastern and western hemispheres
are stable under the $SL(2,\BB R)$-action. So, let $H$ denote one of these
two open hemispheres, say, the western one; and let $S^1 \subset \BB CP^1$
denote the circle containing the Greenwich meridian; $S^1 = \partial H$.

\begin{figure}
\centerline{\psfig{figure=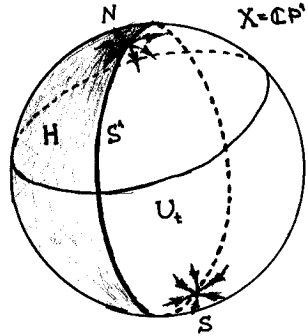}}
\caption{Positioning of $N$, $S$, $H$, $S^1$ and $U_t$.}
\end{figure}

There are only two essential choices for the
$SL(2,\BB R)$-equivariant sheaf ${\cal F}$:
${\cal F}_{\text{principal}}$ -- the sheaf which is
``zero outside of $S^1$ and the constant sheaf along $S^1$'' --
defined by
$$
{\cal F}_{\text{principal}} = (j_{S^1 \hookrightarrow \BB CP^1})_* \BB C_{S^1},
$$
where $\BB C_{S^1}$ denotes the constant sheaf on $S^1$;
and ${\cal F}_{\text{discrete}}$ -- the sheaf which is
``zero outside of $H$ and the constant sheaf along $H$'' --
defined by
$$
{\cal F}_{\text{discrete}} = (j_{H \hookrightarrow \BB CP^1})_! \BB C_H,
$$
where $\BB C_H$ denotes the constant sheaf on $H$.
If we take any other $SL(2,\BB R)$-equivariant sheaf ${\cal F}$,
its characteristic cycle $Ch({\cal F})$ will be an integral linear combination
of $Ch({\cal F}_{\text{principal}})$, $Ch({\cal F}_{\text{discrete}})$
and $Ch({\cal F}_{\text{discrete'}})$, where
$$
{\cal F}_{\text{discrete'}} = (j_{H' \hookrightarrow \BB CP^1})_! \BB C_{H'}
$$
and $H'$ denotes the open hemisphere opposite to $H$ (the eastern hemisphere).

In the case of ${\cal F}_{\text{principal}}$, the integral character formula
(\ref{intformula}) calculates the character of the principal series
representation.
On the other hand,
in the case of ${\cal F}_{\text{discrete}}$, for certain values of $\lambda$,
the integral character formula (\ref{intformula}) calculates the character
of the discrete  series representation.

If $M$ is a submanifold of $X$, we will denote by $T^*_MX$ the
{\em conormal vector bundle} to $M$ in $T^*X$.
That is $T^*_MX$ is a subbundle of $T^*X$ restricted to $M$ consisting of those
covectors in $T^*_mX$ which annihilate $T_m M \subset T_m X$, $m \in M$.

\begin{figure}
\centerline{\psfig{figure=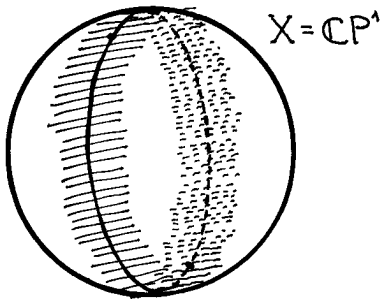}}
\caption{$Ch({\cal F}_{\text{principal}})$.}
\end{figure}

As promised, we give an explicit description of the characteristic cycles:
$$
Ch({\cal F}_{\text{principal}}) = T^*_{S^1} \BB CP^1,
$$
the conormal bundle to $S^1$; and the open embedding theorem of
W.~Schmid and K.~Vilonen (\cite{SchV1}) tells us that
\begin{multline*}
Ch({\cal F}_{\text{discrete}}) = T^*_H \BB CP^1 \bigcup
\{ \zeta \in T^*_{S^1} \BB CP^1;\, \langle \zeta, v \rangle \ge 0  \\
\text{ for all tangent vectors $v$ pointing outside of $H$} \}.
\end{multline*}
Of course, their orientations must be specified too, but they are not
important at this level of details.

\begin{figure}
\centerline{\psfig{figure=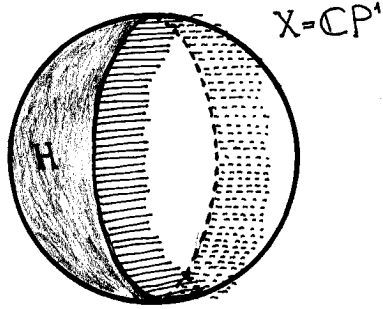}}
\caption{$Ch({\cal F}_{\text{discrete}})$.}
\end{figure}

For $t \in (-90^{\circ},90^{\circ}]$ we will denote by $U_t$ the open
set consisting of all points on the sphere $\BB CP^1$ whose
latitude is strictly less than $t$.

Now we will describe the deformations of $Ch({\cal F}_{\text{discrete}})$
first and then of $Ch({\cal F}_{\text{principal}})$.
For $t \in (-90^{\circ},90^{\circ}]$ we can consider the sheaf
${\cal F}_{\text{discrete}, t}$ which is
``zero outside of $H \cap U_t$ and the constant sheaf along $H \cap U_t$''
defined by
$$
{\cal F}_{\text{discrete}, t} =
(j_{U_t \hookrightarrow \BB CP^1})_! ({\cal F}_{\text{discrete}} |_{U_t})
= (j_{H \cap U_t \hookrightarrow \BB CP^1})_! \BB C_{H \cap U_t}.
$$
Observe that
$$
{\cal F}_{\text{discrete}, 90^{\circ}} = {\cal F}_{\text{discrete}}.
$$
Heuristically, if one ``deforms a sheaf slightly'' (whatever that means)
one could expect its characteristic cycle make only slight changes too.

Then we can consider the characteristic cycle of
${\cal F}_{\text{discrete}, t}$, $Ch({\cal F}_{\text{discrete}, t})$.
If $t=90^{\circ}$,
$Ch({\cal F}_{\text{discrete}, t}) = Ch({\cal F}_{\text{discrete}})$.
When $t \in (-90^{\circ},90^{\circ})$,
we use the open embedding theorem of W.~Schmid and K.~Vilonen (\cite{SchV1})
once again to see that
\begin{multline*}
Ch({\cal F}_{\text{discrete}, t})
= Ch \bigl(
(j_{H \cap U_t \hookrightarrow \BB CP^1})_! \BB C_{H \cap U_t} \bigr)  \\
= T^*_{H \cap U_t} \BB CP^1 \bigcup
\{ \zeta \in T^*_x \BB CP^1;\, x \in \partial (H \cap U_t), \,
\langle \zeta, v \rangle \ge 0  \\
\text{ for all tangent vectors $v$ pointing outside of $H \cap U_t$} \}.
\end{multline*}

\begin{figure}
\centerline{\psfig{figure=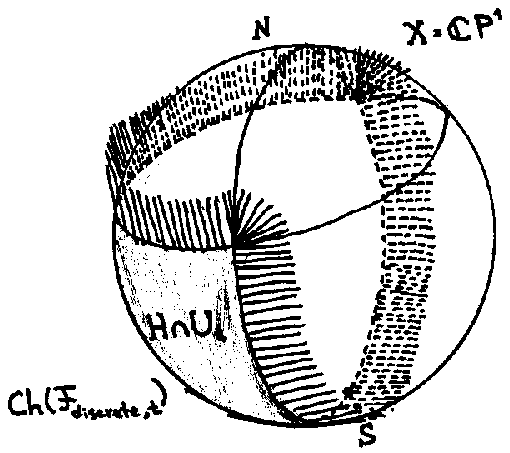}}
\caption{$Ch ({\cal F}_{\text{discrete}, t})$.}
\end{figure}

\begin{figure}
\centerline{\psfig{figure=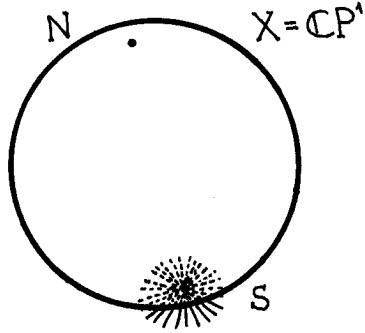}}
\caption{$Ch ({\cal F}_{\text{discrete}})$ is deformed into $T^*_S \BB CP^1$.}
\end{figure}

We see that, as we decrease $t$ from $90^{\circ}$ to $0^{\circ}$,
we continuously deform the cycle $Ch({\cal F}_{\text{discrete}})$.
The end result of this deformation
is the cotangent space at the south pole $S$, $T^*_S \BB CP^1$.
Moreover, because the vector field $\operatorname{VF}_g$ never points
outside of $H \cap U_t$, it is true that during the
process of deformation we always stayed inside the set
$$
\{\zeta \in T^* \BB CP^1;\, Re( \langle g,\mu(\zeta) \rangle) =
\langle \operatorname{VF}_g, \zeta \rangle \le 0\}.
$$
In other words, $Ch({\cal F}_{\text{discrete}}) - T^*_S \BB CP^1$ is the
boundary of a chain which lies entirely inside
\begin{equation}  \label{negative}
\{\zeta \in T^* \BB CP^1;\, Re( \langle g,\mu(\zeta) \rangle) \le 0\}.
\end{equation}

Notice also that although the north pole $N$ and the south pole $S$ apriori
seem to be entirely symmetric with respect to ${\cal F}$, it is only the
south pole $S$ that counts because the condition (\ref{negative}) must be
satisfied during a deformation. And there is no way to deform
$Ch({\cal F}_{\text{discrete}})$ into $T^*_N \BB CP^1$ without breaking
the condition (\ref{negative}).

\separate

Next we show how to deform $Ch({\cal F}_{\text{principal}})$.
First of all, we observe that the flag variety
$$
X = \BB CP^1 = U_{90^{\circ}} \coprod \{N\}.
$$
(For general real reductive Lie groups $G_{\BB R}$ we use
the Bruhat cell decomposition.)
Hence we get a distinguished triangle:
$$
(j_{U_{90^{\circ}} \hookrightarrow \BB CP^1})_!
({\cal F}_{\text{principal}}|_{U_{90^{\circ}}})
\to {\cal F}_{\text{principal}} \to
(j_{\{N\} \hookrightarrow \BB CP^1})_!
({\cal F}_{\text{principal}}|_{\{N\}}),
$$
which is equivalent to
$$
(j_{S^1 \cap U_{90^{\circ}} \hookrightarrow \BB CP^1})_!
\BB C_{U_{S^1 \cap U_{90^{\circ}}}}
\to {\cal F}_{\text{principal}} \to
(j_{\{N\} \hookrightarrow \BB CP^1})_! \BB C_N.
$$
The first term is the sheaf which is ``zero outside
$S^1 \cap U_{90^{\circ}} = S^1 \setminus \{N\}$ and the constant sheaf
along $S^1 \cap U_{90^{\circ}}$''; the last term is the sheaf which is
``zero outside $\{N\}$ and the constant sheaf on $\{N\}$''.

It is a basic property of characteristic cycles that the characteristic
cycle of the middle sheaf in a distinguished triangle equals the sum of the
characteristic cycles of the other two sheaves:
$$
Ch({\cal F}_{\text{principal}}) =
Ch \bigl( (j_{S^1 \cap U_{90^{\circ}} \hookrightarrow \BB CP^1})_!
\BB C_{S^1 \cap U_{90^{\circ}}} \bigr) +
Ch \bigl( (j_{\{N\} \hookrightarrow \BB CP^1})_! \BB C_N \bigr).
$$
The term
$$
Ch \bigl( (j_{\{N\} \hookrightarrow \BB CP^1})_! \BB C_N \bigr) =
T^*_N \BB CP^1.
$$
already has the desired form. We deform the other term
$$
Ch \bigl( (j_{S^1 \cap U_{90^{\circ}} \hookrightarrow \BB CP^1})_!
\BB C_{S^1 \cap U_{90^{\circ}}} \bigr) = T^*_{S^1} \BB CP^1 - T^*_N \BB CP^1
$$
only, and we do it very similarly to the case of
$Ch({\cal F}_{\text{discrete}})$.

Thus, for $t \in (-90^{\circ},90^{\circ}]$, we consider the sheaf
${\cal F}_{\text{principal}, t}$ which is ``zero outside
$S^1 \cap U_t = S^1 \setminus \{N\}$ and the constant sheaf
along $S^1 \cap U_t$'' defined by
$$
{\cal F}_{\text{principal}, t} =
(j_{U_t \hookrightarrow \BB CP^1})_! ({\cal F}_{\text{principal}} |_{U_t})
= (j_{S^1 \cap U_t \hookrightarrow \BB CP^1})_! \BB C_{S^1 \cap U_t}.
$$
Then we consider its characteristic cycle
$Ch({\cal F}_{\text{principal}, t})$.
When $t=90^{\circ}$,
$Ch({\cal F}_{\text{principal}, t}) = T^*_{S^1} \BB CP^1 - T^*_N \BB CP^1$.
When $t \in (-90^{\circ},90^{\circ})$,
we use the open embedding theorem of W.~Schmid and K.~Vilonen (\cite{SchV1})
once again to see that
\begin{multline*}
Ch({\cal F}_{\text{principal}, t})
= Ch((j_{S^1 \cap U_t \hookrightarrow \BB CP^1})_! \BB C_{S^1 \cap U_t})  \\
= T^*_{S^1 \cap U_t} \BB CP^1 \bigcup
\{ \zeta \in T^*_x \BB CP^1;\, x \in \partial (S^1 \cap U_t), \,
\langle \zeta, v \rangle \ge 0  \\
\text{ for all tangent vectors $v \in T_x S^1$ pointing outside of $U_t$} \}.
\end{multline*}

\begin{figure}
\centerline{\psfig{figure=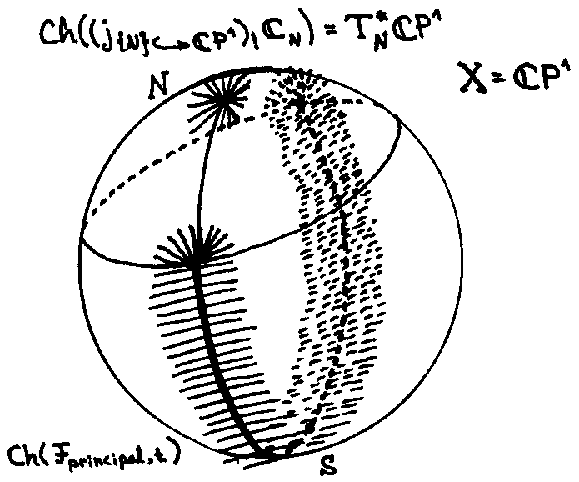}}
\caption{
$Ch \bigl( (j_{\{N\} \hookrightarrow \BB CP^1})_! \BB C_N \bigr)$
and $Ch({\cal F}_{\text{principal}, t})$.}
\end{figure}

\begin{figure}
\centerline{\psfig{figure=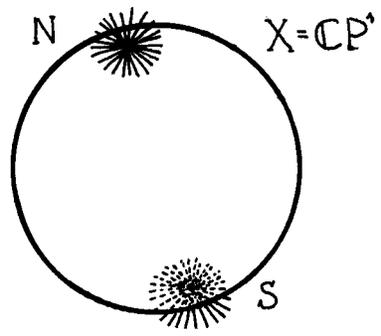}}
\caption{$Ch ({\cal F}_{\text{principal}})$ is deformed into
$T^*_N \BB CP^1 + T^*_S \BB CP^1$.}
\end{figure}

We see that, as we decrease $t$ from $90^{\circ}$ to $0^{\circ}$,
we continuously deform the cycle
$$
Ch \bigl( (j_{S^1 \cap U_{90^{\circ}} \hookrightarrow \BB CP^1})_!
\BB C_{S^1 \cap U_{90^{\circ}}} \bigr).
$$
The end result of this deformation
is the cotangent space at the south pole $S$, $T^*_S \BB CP^1$.
Moreover, as before, it is true that during the
process of deformation we always stayed inside the set
$$
\{\zeta \in T^* \BB CP^1;\, Re( \langle g,\mu(\zeta) \rangle) =
\langle \operatorname{VF}_g, \zeta \rangle \le 0\}.
$$
In other words, $Ch({\cal F}_{\text{principal}}) -
(T^*_N \BB CP^1 + T^*_S \BB CP^1)$ is the
boundary of a chain which lies entirely inside
$$
\{\zeta \in T^* \BB CP^1;\, Re( \langle g,\mu(\zeta) \rangle) \le 0\}.
$$

\separate

Let me remind once again what we have achieved so far.
We started with a cycle $\g{sl}(2,\BB R)' \times Ch({\cal F})$
inside $\g{sl}(2,\BB R)' \times T^* \BB CP^1$ and deformed it
into a new cycle. Let us call it $C_{\cal F}$.
The new cycle $C_{\cal F}$ is characterized by the following two properties:
First of all, for each regular semisimple $g \in \g{sl}(2,\BB R)'$,
$C_{\cal F}$ intersects $\{ g \} \times T^*\BB CP^1$ transversally and
$$
C_{\cal F} \cap (\{ g \} \times T^*\BB CP^1)
= m_1 T^*_{x_1} \BB CP^1 + m_2 T^*_{x_2} \BB CP^1,
$$
where $m_1,m_2$ are some integers,
$x_1,x_2$ are the points in $\BB CP^1$ fixed by $g$,
and each cotangent space $T^*_{x_k} \BB CP^1$ is given some orientation.
Secondly, the deformation of the cycle
$\g{sl}(2,\BB R)' \times Ch({\cal F})$ into $C_{\cal F}$ takes place
inside the set
$$
\{(g,\zeta) \in \g{sl}(2,\BB R)' \times T^* \BB CP^1;\,
Re( \langle g,\mu(\zeta) \rangle) \le 0\}.
$$
These two properties determine the cycle $C_{\cal F}$ uniquely.

\separate

The way this deformation procedure is done in general, it becomes
very similar to the classical Morse's lemma which says
that if we have a smooth real valued function $f$ on a manifold
$M$, then the sublevel sets $\{ m \in M;\,f(m) < a \}$ and
$\{ m \in M;\,f(m) < b \}$ can be deformed one into the other
as long as there are no critical values of $f$ in an open
interval containing $a$ and $b$.

\separate

After deforming the cycle
$\g{sl}(2,\BB R)' \times Ch({\cal F})$ into $C_{\cal F}$
we encounter the following obstacle: The integrand
$e^{\langle g, \mu_{\lambda}(\zeta) \rangle}
\phi(-\sigma+\pi^* \tau_{\lambda})$
when restricted to $C_{\cal F}$ becomes zero!
There is no contradiction to Stokes' theorem because
the cycles in question are not compact.
In this situation we are simply not allowed to interchange pointwise limits
and integration. See Remark \ref{zero_remark}.

In order to deal with this problem, for each regular semisimple
$g \in \g{sl}(2,\BB R)'$, we define another deformation
$\Theta_t(g): T^*X \to T^*X$, $t \in [0,1]$.
Its description will be given right after the equation (\ref{RR}).
These $\Theta_t(g)$'s combine into a deformation
$$
\Theta_t: \g{sl}(2,\BB R)' \times T^*X \to \g{sl}(2,\BB R)' \times T^*X,
\qquad t \in [0,1].
$$
The idea of this deformation was inspired by the classical proof of
the Fourier inversion formula
$$
\phi(g)= \frac 1{(2\pi i)^{\dim_{\BB C}\g g}}
\int_{\zeta \in i \g g_{\BB R}^*}
\hat \phi(\zeta) e^{-\langle g,\zeta \rangle}
$$
where we multiply the integrand by a term like $e^{-t \|\zeta\|^2}$
to make it integrable over $\g g_{\BB R} \times i \g g_{\BB R}^*$,
and then let $t \to 0^+$.
The deformation $\Theta_t$ has a very similar effect --
it makes our integrand an $L^1$-object.
The following result is Lemma 17 in \cite{L}:

\begin{lm} \label{slanting}
For any $t \in [0,1]$, we have:
\begin{multline*}
\lim_{R \to \infty}
\int_{\g{sl}(2,\BB R)' \times (Ch({\cal F}) \cap \{\|\zeta\| \le R\})}
\bigl( e^{\langle g, \mu_{\lambda}(\zeta) \rangle} \phi
(-\sigma + \pi^* \tau_{\lambda})   \\
- \Theta_t^*(e^{\langle g, \mu_{\lambda}(\zeta) \rangle} \phi
(-\sigma + \pi^* \tau_{\lambda})) \bigr) =0.
\end{multline*}
\end{lm}

It essentially says that it is permissible to substitute
$$
\Theta_t^*(e^{\langle g, \mu_{\lambda}(\zeta) \rangle} \phi
(-\sigma + \pi^* \tau_{\lambda}))
$$
in place of the original integrand
$e^{\langle g, \mu_{\lambda}(\zeta) \rangle} \phi
(-\sigma + \pi^* \tau_{\lambda})$.

Proof of Lemma \ref{slanting} is very technical, but the idea is quite simple.
The difference between the original integral and the deformed
one is expressed by an integral of 
$e^{\langle g, \mu_{\lambda}(\zeta) \rangle}
\phi(-\sigma+\pi^* \tau_{\lambda})$ over a certain cycle
$\tilde C(R)$ supported in
$\g{sl}(2,\BB R)' \times (T^*\BB CP^1 \cap \{\|\zeta\| = R\})$
which depends on $R$ by scaling along the fiber.
Recall that the Fourier transform $\hat \phi$ decays rapidly
in the imaginary directions which is shown by an integration
by parts. We modify this integration by parts argument to prove
a similar statement about behavior of the integrand on
the support of $\tilde C(R)$ as $R \to \infty$.
Hence the difference of integrals in question tends to zero.

\separate

Let $C$ be the Borel-Moore chain in $\g{sl}(2,\BB R)' \times T^*\BB CP^1$
of dimension 6 which lies inside the set
$$
\{(g,\zeta) \in \g{sl}(2,\BB R)' \times T^* \BB CP^1;\,
Re( \langle g,\mu(\zeta) \rangle) \le 0\}
$$
and such that
$$
\partial C = \g{sl}(2,\BB R)' \times Ch({\cal F}) - C_{\cal F}.
$$
Take an $R \ge 1$ and restrict all cycles to the set
$$
\{ (g,\zeta) \in \g{sl}(2,\BB R)' \times T^*\BB CP^1;\, \|\zeta\| \le R\}.
$$
Then the restricted chain $C$ has boundary
\begin{multline*}
\partial (C \cap \{\|\zeta\| \le R\}) =
\g{sl}(2,\BB R)' \times (Ch({\cal F}) \cap \{\|\zeta\| \le R\})  \\
- (C_{\cal F} \cap \{\|\zeta\| \le R\}) - C'(R),
\end{multline*}
where $C'(R)$ is a 5-chain supported in the set
$$
\{ (g, \zeta) \in \g{sl}(2,\BB R)' \times T^*X;\, \|\zeta\|=R,\,
Re( \langle g, \mu(\zeta) \rangle ) \le 0 \}.
$$
Because our chain $C$ is conic, the piece of boundary $C'(R)$
depends on $R$ by appropriate scaling of $C'(1)$ in the fiber direction.

The following result is Lemma 18 in \cite{L}:

\begin{lm}  \label{zero}
For a fixed $t \in (0,1]$,
$$
\lim_{R \to \infty} \int_{C'(R)} 
\Theta_t^* \bigl( e^{\langle g, \mu_{\lambda}(\zeta) \rangle} \phi
(-\sigma + \pi^* \tau_{\lambda})^n \bigr) =0.
$$
\end{lm}

\separate

We start with the right hand side of the
integral character formula (\ref{intformula}) and
continuing calculations (\ref{R}) using Lemmas \ref{slanting} and \ref{zero}
we obtain:
\begin{multline}  \label{RR}
\frac 1{2\pi i} \int_{Ch({\cal F})}
\mu_{\lambda}^* \hat \phi (-\sigma+\pi^* \tau_{\lambda})  \\
=\lim_{R \to \infty} \frac 1{2\pi i}
\int_{\g{sl}(2,\BB R)' \times (Ch({\cal F}) \cap \{\|\zeta\| \le R\})}
e^{\langle g, \mu_{\lambda}(\zeta) \rangle} \phi
(-\sigma+\pi^* \tau_{\lambda})  \\
=\lim_{R \to \infty} \frac 1{2\pi i}
\int_{\g{sl}(2,\BB R)' \times (Ch({\cal F}) \cap \{\|\zeta\| \le R\})}
\Theta_t^*(e^{\langle g, \mu_{\lambda}(\zeta) \rangle} \phi
(-\sigma + \pi^* \tau_{\lambda}))  \\
=\lim_{R \to \infty} \frac 1{2\pi i}
\int_{(C_{\cal F} \cap \{\|\zeta\| \le R\}) + C'(R)}
\Theta_t^*(e^{\langle g, \mu_{\lambda}(\zeta) \rangle} \phi
(-\sigma + \pi^* \tau_{\lambda}))  \\
=\lim_{R \to \infty} \frac 1{2\pi i}
\int_{C_{\cal F} \cap \{\|\zeta\| \le R\}}
\Theta_t^*(e^{\langle g, \mu_{\lambda}(\zeta) \rangle} \phi
(-\sigma + \pi^* \tau_{\lambda}))  \\
= \frac 1{2\pi i} \int_{C_{\cal F}}
\Theta_t^* \bigl( e^{\langle g, \mu_{\lambda}(\zeta) \rangle} \phi
(-\sigma + \pi^* \tau_{\lambda}) \bigr).
\end{multline}
It turns out that the last integral is pretty easy to calculate.

\separate
As before,
$g \in \g{sl}(2,\BB R)'$ is a regular semisimple element, and
$x_1,x_2 \in \BB CP^1$ are the zeroes of the vector field
$\operatorname{VF}_g$ on $\BB CP^1$.
For $k =1,2$, we let $\g b_k \subset \g g$ be the Borel subalgebra
of $\g g$ consisting of all elements in $\g g$ fixing $x_k$.
Let $\g n_k \subset \g g$ be the root
spaces with respect to the unique Cartan algebra $\g t \subset \g g$ 
containing $g$ such that $\g g = \g b_k \oplus \g n_k$ as linear spaces.
Define maps $\psi_k: \g n_k \to \BB CP^1$, $n \mapsto \exp(n) \cdot x_k$.
Each $\psi_k$ is a diffeomorphism onto its image and their images
cover all of $\BB CP^1$, i.e. $\{\psi_1,\psi_2\}$ is an atlas of $\BB CP^1$.
Let $z_k: \g n_k \tilde \to \BB C$ be a linear coordinate on $\g n_k$.
Then it is not hard to show that the vector field $\operatorname{VF}_g$
can be expressed in these coordinates by
$$
\alpha_{x_k}(g) z_k \frac{\partial}{\partial z_k},
$$
where $\alpha_k \in \g t^*$ denotes the root of $\g{sl}(2,\BB C)$
corresponding to the root space $\g n_k$.

We expand $z_k$ to a standard coordinate system $z_k, \xi_k$ on
$T^* \g n_k$ so that every element of
$T^* \g n_k \simeq \g n_k \times \g n_k^*$ is expressed in these
coordinates as $(z_k, \xi_k dz_k)$. This gives us two charts
$\tilde \psi_k : (z_k,\xi_k) \to T^*\BB CP^1$ and an atlas
$\{\tilde \psi_1, \tilde \psi_2 \}$ of $T^*\BB CP^1$.

The diffeomorphism $\Theta_t(g): T^*\BB CP^1 \to T^*\BB CP^1$ is defined
using the atlas $\{\tilde \psi_1, \tilde \psi_2 \}$ and smooth cutoff
functions. It is designed so that it approximates the maps
\begin{equation}  \label{Theta}
(z_k, \xi_k) \mapsto
\bigl( z_k - t \frac{\bar \alpha_k(g)}{|\alpha_k(g)|} \bar \xi_k,\xi_k \bigr).
\end{equation}
Of course, these maps cannot be defined globally on $T^*\BB CP^1$, but for
the purposes of calculating the integral (\ref{RR}) we are allowed
to pretend that, in each chart $\tilde \psi_k$, $\Theta_t(g)$ is given
by the equation (\ref{Theta}).

Next we calculate the restriction of the integrand of (\ref{RR})
$$
\Theta_t(g)^* \bigl( e^{\langle g, \mu_{\lambda}(\zeta) \rangle} \phi
(-\sigma + \pi^* \tau_{\lambda}) \bigr)
$$
to $T^*_{x_k} \BB CP^1$:
$$ 
\Theta_t(g)^*
\bigl( \langle g, \mu_{\lambda}(\zeta) \rangle \bigr) |_{T^*_{x_k} \BB CP^1}
= - t |\alpha_k(g)| \cdot |\xi_k|^2
+ \Theta_t(g)^* \langle g, \lambda_x \rangle,
$$
$$
\Theta_t(g)^* \sigma |_{T^*_{x_k} \BB CP^1}
= t \frac{\bar \alpha_k(g)}{|\alpha_k(g)|} d\bar \xi_k \wedge d\xi_k,
$$
and $\Theta_t(g)^* (\pi^* \tau_{\lambda})$ on $T^*_{x_k} \BB CP^1$ 
can be expressed as $t^2$ times a bounded $2$-form.
Putting it all together, we obtain
\begin{multline*}
\Theta_t(g)^* \bigl( e^{\langle g, \mu_{\lambda}(\zeta) \rangle} \phi
(-\sigma + \pi^* \tau_{\lambda}) \bigr) |_{T^*_{x_k} \BB CP^1}   \\
= e^{- t |\alpha_k(g)| \cdot |\xi_k|^2
+ \Theta_t(g)^* \langle g, \lambda_x \rangle}
\phi(g) \wedge
\Bigl( t \frac{\bar \alpha_k(g)}{|\alpha_k(g)|} d\xi_k \wedge d\bar\xi_k
+ t^2 \cdot (\text{bounded term}) \Bigr).
\end{multline*}
One can argue using a coordinate change and the Lebesgue dominated
convergence theorem that, as $t \to 0^+$, we can replace the integrand with
\begin{equation}  \label{integrand_t}
t \frac{\bar \alpha_k(g)}{|\alpha_k(g)|}
e^{\langle g, \lambda_{x_k} \rangle} 
e^{- t |\alpha_k(g)| \cdot |\xi_k|^2}
\phi(g) \wedge d\xi_k \wedge d\bar\xi_k,
\end{equation}
which integrates over $T^*_{x_k} \BB CP^1$ to
$$
(2\pi i) \frac {e^{\langle g, \lambda_{x_k} \rangle}}{\alpha_k(g)} \phi(g).
$$
The last expression may appear to be missing a negative sign, but it is
correct because we are using the convention (11.1.2) in \cite{KaScha}
to identify the real cotangent bundle $T^*(X^{\BB R})$ with the
holomorphic cotangent bundle $T^*X$ of the complex manifold $X$.
Under this identification, the standard symplectic
structure on $T^*(X^{\BB R})$ equals $2 Re (\sigma)$.

Substituting into (\ref{RR}) we obtain
$$
\frac 1{2\pi i} \int_{Ch({\cal F})}
\mu_{\lambda}^* \hat \phi (-\sigma+\pi^* \tau_{\lambda})
= \int_{\g g_{\BB R}} \sum_{k=1}^2 m_k(g)
\frac {e^{\langle g, \lambda_{x_k(g)} \rangle}}{\alpha_{k(g)}(g)} \phi(g)
= \int_{\g g_{\BB R}} F_{\theta} \phi,
$$
which says that the right hand side of the integral character formula
(\ref{intformula}) is equal to
the right hand side of the fixed point character formula (\ref{fpformula}).

\begin{rem} \label{zero_remark} {\em
Notice that, as $t \to 0^+$, the integrand (\ref{integrand_t}) tends to zero
pointwise, but its integral does not.
We had to introduce the deformation $\Theta_t(g) : T^* X \to T^*X$ so that
we could interchange pointwise limits and integrals and resolve the
situation when integrand (\ref{integrand}) restricted to the cycle
$C_{\cal F}$ is zero, but the integral does not tend to zero.
}\end{rem}

\separate

The key ideas are the deformation of $Ch({\cal F})$,
the definition of $\Theta_t(g): T^*X \to T^*X$ and
Lemma \ref{slanting}. Because of the right definition
of $\Theta_t(g)$, Lemma \ref{slanting} holds and our calculation
of the integral becomes very simple.
We see that, as $R \to \infty$,
the integral concentrates more and more inside $T^*U$,
where $U$ is a neighborhood of the set of fixed points of $g$ in $X$.
In the limit, we obtain the right hand side of the
fixed point character formula.
This means that the integral is localized at the fixed points of $g$.

\separate

As was stated before, complete proofs with all technical details
are available in \cite{L}.

\separate


\noindent
{matvei@math.umass.edu}

\noindent
{\em Department of Mathematics and Statistics, University of Massachusetts,
Lederle Graduate Research Tower, 710 North Pleasant Street, Amherst,
MA 01003}

\enddocument